\documentclass[a4paper,12pt]{article}
\usepackage{amsmath,amsfonts,amssymb,latexsym,epsfig,mathrsfs}
\usepackage[all]{xy}
\usepackage{graphicx}
\newcommand{\be}{\begin{equation}}
\newcommand{\ee}{\end{equation}}

\newtheorem{theorem}{Theorem}[subsection]
\newtheorem{definition}{Definition}[subsection]
\newtheorem{proposition}{Proposition}[subsection]

\newtheorem{remark}{Remark}[subsection]

\def\proof{\bf Proof. \rm}

\begin{document}
\thispagestyle{empty}
 \begin{center}
{{\bf INDUCED RIEMANNIAN STRUCTURE ON \\A REDUCED SYMPLECTIC MANIFOLD}\\

\vspace{0.5cm}
	
\small {Augustin T. Batubenge* \quad\quad\quad\quad\quad \quad\quad\quad Wallace M. Haziyu  }}
\end{center}

\vspace{1cm}

\begin{abstract}

Let $(M,\omega)$ be a symplectic and Riemannian manifold on which a Lie group acts in a hamiltonian way. We determine conditions for which the reduced space inherits an induced Riemannian structure through a Riemannian submersion.
\end{abstract}

\section{Introduction}
These notes aim at studying the inheritance of a Riemann metric of a symplectic manifold on its symplectic quotient. In several studies the symplectic quotient process has been investigated throughout the last century, which aimed at working out a lower dimensional manifold that still holds a symplectic structure related to the original one. In mechanics the new object offers that advantage of modeling a Hamiltonian system with lower degree of freedom. Our main references on the topic are the paper by Marsden G. and Weinstein A. (\cite[1974]{Mar74}) and the book by Abraham R. (\cite{Abr78}). In the same way, starting with a symplectic manifold provided with a Rieman metric, we think it would be good to end up having a Marsden-Weinstein quotient which is also a Riemannian space with a Riemann metric inherited from the one on the original space. To this end, the paper is organized as follows.\\
Section 1 is a summary of the basics on Riemannian manifolds that introduces Riemannian submersions. Section 2 describes almost complex structures and the notion of a holomorphic mapping that will play a key role in the transference of the metric to the quotient, that is, the existence of a Riemann structure on the symplectic quotient through a Riemann submersion.      

\section{Preliminary}

A Riemannian structure (or Riemannian metric) on a smooth manifold $M$, which is usually denoted by $g$, is a smooth positive definite and symmetric bilinear form on $T_pM$ for each $p\in M$. That is, a smooth assignment of an inner product $\langle \cdot,\cdot\rangle$, to each tangent space $T_{p}M$ of $M$. We denote by $(M,g)$ a manifold on which the Riemannian structure $g$ is defined and call it the Riemannian manifold. The inner product at each $p\in M$ we shall denote by either $g_p$ or $g_p(M)$. 

\vspace{0.3cm}

Assuming that $N$ is an arbitrary manifold, $(M,g)$ a Riemannian manifold, we recall that if a smooth map $f:N\rightarrow M$ is an immersion, then $f^{*}g$ is a Riemannian metric on $N$ called the induced metric. Moreover, let $(M, g)$ and $(N , h)$ be two Riemannian manifolds. A diffeomorphism $f:M\rightarrow N$ is called an isometry if
$g_{p}(X , Y)$ = $h_{f(p)}(T_{p}f\cdot X , T_{p}f\cdot Y)$ for all $X,Y\in T_{p}M$, $p\in M$ and where $T_{p}f\cdot X$ is the image of the tangent vector $X$ by the differential mapping associated with $f$ at $p$. We also say that $f:M\rightarrow M$ is an isometry on $M$ if for all $u , v\in T_{p}M$, $p\in M$, we have $g_{p}(u , v) = g_{f(p)}(T_{p}f\cdot u , T_{p}f\cdot v)$.\\

It is easily checked that if $f$ is an isometry on $M$, then its inverse is also an isometry on $M$. Clearly the identity map on $M$ is an isometry on $M$ and if $f, g$ are isometries on $M$ then their composition is also an isometry on $M$, which makes the set of isometries on $M$ a group under the composition of maps.
A group of isometries on a Riemannian manifold $M$ is a Lie group. (see \cite [p~63] {G-H-L87})

\begin{theorem}\label{th:th01}
Let $G$ be a Lie group of isometries of a Riemannian manifold $(M, g)$ acting transitively on $M$, then $G$ is compact if and only if $M$ is compact.

\end{theorem}

For the proof of this theorem (see \cite [Theorem~2.35, p.~63]{G-H-L87}).

\begin{definition}

Let $\Phi:G\times M\rightarrow M$ be an action of a Lie group $G$ on a smooth manifold $M$. Then a Riemannian metric $g(\cdot , \cdot)$ on $M$ is called invariant if for each $m\in M$ we have $g_{m}(u , v) = g_{\Phi_{a}(m)}(T_{m}\Phi_{a}\cdot u , T_{m}\Phi_{a}\cdot v)$ for all $u , v\in T_{m}M$ and $a\in G$.
\end{definition}

\begin{theorem}\label{th:th02}

Let $G$ be a Lie group acting on a smooth manifold $M$. If $G$ is compact then there exists an invariant Riemannian metric on $M$ (see \cite [p.~56] {Aud04}).
\end{theorem}

\vspace{0.3cm}

\subsection{Riemannian submersions}

\begin{definition}
Let $f:M\rightarrow N$ be a smooth map. An element $x\in N$ is called a regular value of $f$ if $f^{-1}(x)$ is a submanifold of $M$, and if whenever $m\in f^{-1}(x)$ then $T_{m}f:T_{m}M\rightarrow T_{f(m)}N$ is surjective. A point $m\in M$ is called a regular point of $f$ if $T_{m}f$ is surjective.
\end{definition}

\begin{definition}

Let $M$ and $N$ be smooth manifolds. A smooth map $\Phi:M\rightarrow N$ is called a submersion if all points of $M$ are regular points of $\Phi$. That is, $\Phi$ is a submersion if
\[
\begin{array}{ccc}
(d\Phi)_{x}: T_{x}M\rightarrow T_{\Phi(x)}N
\end{array}
\]
is surjective for all $x\in M$
\end{definition}

\begin{remark}
If $\Phi:M\rightarrow N$ is a submersion then each $b\in N$ is a regular value of $\Phi$ and $\Phi^{-1}(b)$ is called a fibre of $\Phi$ over {b}. The rank of $\Phi$ is equal to the dimension of $N$, i.e rank$\Phi = \dim N$ at every point of $\Phi^{-1}(b)$ for all $b\in N$. It is known that $\Phi^{-1}(b)$ is a regular submanifold of $M$ of dimension $\dim M-\dim N$ (\cite [p.~391]{Che05}) and (\cite [p.~459] {Bar66}). If $p\in \Phi^{-1}(b)$ then $T_{p}\Phi^{-1}(b)$ is the kernel of the differential of $\Phi$ at $p$. That is, $T_{p}\Phi^{-1}(b) = \ker{d\Phi}_{p}$. So if $v_{p}\in T_{p}\Phi^{-1}(b)$ then $d{\Phi}_{p}\cdot v_{p} = 0$. Let $V(M)_{p} = T_{p}\Phi^{-1}(b)$. The set $V(M)_{p}$ is called the set of vertical vectors at $p$.
\end{remark}

\begin{definition}
If $(M , g)$ and $(B , h)$ are Riemannian manifolds, a smooth map
\[
\begin{array}{ccc}
\pi: M\rightarrow B
\end{array}
\]
is called a Riemannian submersion if:

\begin{description}
\item   (i) $\pi$ has maximum rank at each point $p\in M$. That is to say $(d\pi)_{p}:T_{p}M\rightarrow T_{\pi(p)}B$ is surjective, and if we let $V(M)_{p} := \ker (d\pi)_{p}$, then\\
\item   (ii) $(d\pi)_{p}$ is an isometry between $H(M)_{p}$ and $T_{\pi(p)}B$, where $H(M)_{p}$ is the orthogonal complement of $V(M)_{p}$. That is, if $X_{p}, Y_{p}\in H(M)_{p}$ then $g_{p}(X_{p} , Y_{p}) = h_{\pi(p)}((d\pi)_{p}X_{p} , (d\pi)_{p}Y_{p})$.

\end{description}

\end{definition}

We shall denote by $V(M)_{p}$ the set of vertical vectors, and by $H(M)_{p}$ the set of horizontal vectors and note that the tangent space $T_{p}M$ decomposes into an orthogonal direct sum
\[
\begin{array}{cc}
T_{p}M = H(M)_{p}\oplus V(M)_{p},
\end{array}
\]
where $H(M)_{p}\cap V(M)_{p} = \lbrace 0\rbrace$.

\begin{definition}
Let $\Phi:G\times M\rightarrow M$ be the action of a Lie group $G$ on a manifold $M$. The isotropy subgroup of $m\in M$ is a set $G_{m}=\lbrace g\in G:\Phi_{g}(m)= m\rbrace$. The action $\Phi$ is said to be free if whenever $\Phi_{g}(m) = m$ for some $g\in G$ and $m\in M$, then $g = e$, the identity element of $G$. That is, the action is free if the isotropy subgroup is trivial. The action $\Phi$ is called proper if the map $G\times M\rightarrow M\times M$; $\Phi(g,m)\mapsto (\Phi_{g}(m),m)$ is proper. That is, if the inverse image of a compact set is compact.
\end{definition}

\begin{proposition}
Let $G$ be a Lie group of isometries acting properly and freely on a Riemannian manifold $(M , g)$ and let $p:M\rightarrow M/G$ be the canonical projection map (note that N = M/G is a manifold). Then there exists a unique metric on $N = M/G$ such that the projection map $p$ is a Riemannian submersion.
\end{proposition}
(See \cite [p.~61] {G-H-L87}).

\vspace{0.2cm}

We now make the following observations.
\begin{description}

\item   (a) Let $b\in N = M/G$, if $m_{1}, m_{2}\in p^{-1}(b)$ then there is $h\in G$ such that $\Phi_{h}(m_{1}) = m_{2}$ (\cite [Proposition 2.28]{G-H-L87}), where $\Phi$ is the action of $G$ on $M$. Thus the isometry group $G$ acts transitively on each fibre so that the action of $G$ preserves the fibres.

\item   (b) Let $x\in p^{-1}(b)$. For each $\xi\in\mathfrak{g} = T_{e}G$, let $F(t) = \exp{t\xi}$ be its flow, then $\xi_{M}(x) = \frac{d}{dt}\Phi(\exp{t\xi},x)\mid_{t = 0}$ is a tangent vector to the fibre through $x$. If $\xi \neq 0$ then $\xi_{M}(x)\neq 0$. Thus there is a one-to-one correspondence between $\mathfrak{g} = T_{e}G$ and the tangent space to the fibre at each point $x$ in the fibre.

\item   (c) The action of $G$ preserves the vertical distribution $V(M)$. To see this let $\xi_{M}(\cdot)$ be a vertical vector, then for $a\in G$, by straight forward calculations it can be shown that

\[
\begin{array}{cll}

        T_{\Phi_{a^{-1}}(x)}\Phi_{a}\xi_{M}(\Phi_{a^{-1}}(x)) &=& T_{\Phi_{a^{-1}}(x)}\Phi_{a}\frac{d}{dt}\Phi(\exp{t\xi},\Phi_{a^{-1}}(x))\mid_{t = 0}\\

        &=&(Ad_{a}\xi)_{M}(x)
\end{array}
\]
which is a tangent vector to the fibre through $x$.
\end{description}

\vspace{0.2cm}

\subsection{Almost Complex Structure}

Let $\mathbb{C}^{n}$ denote n-dimensional space of complex numbers $(z^{1},z^{2},\cdots , z^{n})$. We identify $\mathbb{C}^{n}$ with $\mathbb{R}^{2n}$ by the correspondence $(z^{1},\cdots ,z^{n})\rightarrow (x^{1},y^{1},\cdots , x^{n},y^{n})$, with $z^{k} = x^{k}+iy^{k}$ where $i = \sqrt{-1}$. By this identification we can consider $\mathbb{C}^{n}$ as a $2n-$dimensional Euclidean space. Similarly, if $M$ is an $n-$ dimensional complex manifold with local coordinates $(z^{1},\cdots , z^{n})$, by identifying these coordinates with $(x^{1},y^{1},\cdots , x^{n},y^{n})$ where $z^{k} = x^{k}+iy^{k}$, $i = \sqrt{-1}$, $k=1,\cdots,n$, we can regard $M$ to be a $2n-$ dimensional differentiable manifold. Then for $p\in M$, the tangent space $T_{p}M$ has the  basis $\lbrace(\frac{\partial}{\partial x^{1}})_{p},(\frac{\partial}{\partial y^{1}})_{p},\cdots, (\frac{\partial}{\partial x^{n}})_{p},(\frac{\partial}{\partial y^{n}})_{p}\rbrace$\\

We now define a linear map $J_{p}:T_{p}M\rightarrow T_{p}M$ by requiring that
\[
\begin{array}{cc}
J_{p}\left(\frac{\partial}{\partial x^{i}}\right)_{p} = \left(\frac{\partial}{\partial y^{i}}\right)_{p}\\

J_{p}\left(\frac{\partial}{\partial y^{i}}\right)_{p} = -\left(\frac{\partial}{\partial x^{i}}\right)_{p}
\end{array}
\]

$i = 1 , 2,\cdots , n$.\\

Clearly $J_{p}^{2} = -\mathbf{1}$.\\

This definition of $J_{p}$ does not depend on the choice of local coordinates $(z^{1},\cdots , z^{n})$. (See \cite [p.~107] {Mat72}).

\begin{definition}

Let $M$ be a smooth manifold, an almost complex structure on $M$ is a correspondence $J$ which assigns a linear transformation\\
$J_{p}:T_{p}M\rightarrow T_{p}M$ for each $p\in M$, with $J^{2}_{p} = -id_{T_{p}M}$, for all $p\in M$.
\end{definition}

The pair $(M , J)$ is called an almost complex manifold.
\vspace{0.3cm}

\begin{proposition}
A differentiable map $\phi:M_{1}\rightarrow M_{2}$ between two almost complex manifolds $M_{1}$ and $M_{2}$ with almost complex structures $J_{1}$ and $J_{2}$ respectively is holomorphic if and only if $\phi_{*}\circ J_{1} = J_{2}\circ \phi_{*}$, where $\phi_{*}$ is the differential of the map $\phi$.
\end{proposition}

Recall that if $F:M\rightarrow N$ is a smooth map and let $\varphi = (x^{1},\cdots , x^{n})$ be local coordinates about $p\in M$ and $\psi = (y^{1},\cdots , y^{m})$ local coordinates about $F(p)\in N$. Then

\[
\begin{array}{ccc}
F_{*}\left(\frac{\partial}{\partial x^{i}}\right)_{p} &=& \displaystyle\sum_{j=1}^{m}\left(F_{*}\left(\frac{\partial}{\partial x^{i}}\right)_{p}y^{j}\right)\frac{\partial}{\partial y^{j}}\mid_{F(p)}\\

&=& \displaystyle\sum_{j=1}^{m}\left(\frac{\partial}{\partial x^{i}}\right)_{p}(y^{j}\circ F)\frac{\partial}{\partial y^{j}}\mid_{F(p)}\\

&=& \displaystyle\sum_{j=1}^{m}\frac{\partial(y^{j}\circ F)}{\partial x^{i}}(p)\frac{\partial}{\partial y^{j}}\mid_{F(p)}
\end{array}
\]

If $f$ is a smooth function on $N$ then the pull back of $f$ under $F$ is a smooth function on $M$ given by $F^{*}f = f\circ F$.\\

\proof (of the Proposition):\\

Let $p\in M_{1}$ and let $(z^{1},\cdots , z^{n})$ be the complex local coordinates in the neighborhood of $p$ and identify these coordinates with $(x^{1},y^{1}, \cdots , x^{n},y^{n})$ of $\mathbb{R}^{2n}$. Let $(w^{1}, \cdots , w^{m})$ be the local coordinates of the neighborhood of $\phi(p)$ in $M_{2}$ identified with $(u^{1},v^{1}, \cdots , u^{m},v^{m})$ of $\mathbb{R}^{2m}$ where
\[
\begin{array}{cll}
z^{k} &=& x^{k}+iy^{k} \hspace{1cm} k = 1, 2, \cdots, n\\
w^{j} &=& u^{j}+iv^{j} \hspace{1cm} j = 1, 2, \cdots, m
\end{array}
\]

Set
\[
\begin{array}{cll}
\phi^{*}u^{j} &=& a_{j}(x^{1},y^{1},\cdots,x^{n},y^{n})\textrm{~and}\\
\phi^{*}v^{j} &=& b_{j}(x^{1},y^{1},\cdots,x^{n},y^{n}), j = 1,\cdots, m
\end{array}
\]
Then by the above comments we have
\[
\begin{array}{cll}
\phi_{*}\displaystyle\left(\frac{\partial}{\partial x^{i}}\right)_{p} &=& \displaystyle\sum_{j=1}^{m}\frac{\partial(u^{j}\circ\phi)}{\partial x^{i}}(p)\frac{\partial}{\partial u^{j}}|_{\phi(p)} + \displaystyle\sum_{j=1}^{m}\frac{\partial(v^{j}\circ\phi)}{\partial x^{i}}(p)\frac{\partial}{\partial v^{j}}|_{\phi(p)}\\

&=& \displaystyle\sum_{j=1}^{m}\frac{\partial a_{j}}{\partial x^{i}}(p)\frac{\partial}{\partial u^{j}}|_{\phi(p)} + \displaystyle\sum_{j=1}^{m}\frac{\partial b_{j}}{\partial x^{i}}(p)\frac{\partial}{\partial v^{j}}|_{\phi(p)}
\end{array}
\]

Similarly
\[
\begin{array}{cll}
\phi_{*}\displaystyle\left(\frac{\partial}{\partial y^{i}}\right)_{p} &=& \displaystyle\sum_{j=1}^{m}\frac{\partial a_{j}}{\partial y^{i}}(p)\frac{\partial}{\partial u^{j}}|_{\phi(p)} + \displaystyle\sum_{j=1}^{m}\frac{\partial b_{j}}{\partial y^{i}}(p)\frac{\partial}{\partial v^{j}}|_{\phi(p)}.
\end{array}
\]

Now from
\begin{equation}\label{eq 01:eq 01}
\phi_{*}\left(J_{1}\frac{\partial}{\partial x^{i}}\right)_{p} = \phi_{*}\left(\frac{\partial}{\partial y^{i}}\right)_{p} = \displaystyle\sum_{j=1}^{m}\frac{\partial a_{j}}{\partial y^{i}}(p)\frac{\partial}{\partial u^{j}}|_{\phi(p)} + \displaystyle\sum_{j=1}^{m}\frac{\partial b_{j}}{\partial y^{i}}(p)\frac{\partial}{\partial v^{j}}|_{\phi(p)} 
\end{equation}

and
\begin{equation}\label{eq 02:eq 02}
\phi_{*}\left(J_{1}\frac{\partial}{\partial y^{i}}\right)_{p} = -\phi_{*}\left(\frac{\partial}{\partial x^{i}}\right)_{p} = -\displaystyle\sum_{j=1}^{m}\frac{\partial a_{j}}{\partial x^{i}}(p)\frac{\partial}{\partial u^{j}}|_{\phi(p)} - \displaystyle\sum_{j=1}^{m}\frac{\partial b_{j}}{\partial x^{i}}(p)\frac{\partial}{\partial v^{j}}|_{\phi(p)}
\end{equation}

On the other hand
\[
\begin{array}{cll}
J_{2}\circ\phi_{*}\displaystyle\left(\frac{\partial}{\partial x^{i}}\right)_{p} = \displaystyle\sum_{j=1}^{m}\frac{\partial a_{j}}{\partial x^{i}}(p)J_{2}\left(\frac{\partial}{\partial u^{j}}\right)_{\phi(p)} + \displaystyle\sum_{j=1}^{m}\frac{\partial b_{j}}{\partial x^{i}}(p)J_{2}\left(\frac{\partial}{\partial v^{j}}\right)_{\phi(p)}
\end{array}
\]
\begin{equation}\label{eq 03:eq 03}
= \displaystyle\sum_{j=1}^{m}\frac{\partial a_{j}}{\partial x^{i}}(p)\frac{\partial}{\partial v^{j}}|_{\phi(p)} - \displaystyle\sum_{j=1}^{m}\frac{\partial b_{j}}{\partial x^{i}}(p)\frac{\partial}{\partial u^{j}}|_{\phi(p)}
\end{equation}

and

\[
\begin{array}{cll}
J_{2}\circ\phi_{*}\displaystyle\left(\frac{\partial}{\partial y^{i}}\right)_{p} = \displaystyle\sum_{j=1}^{m}\frac{\partial a_{j}}{\partial y^{i}}(p)J_{2}\left(\frac{\partial}{\partial u^{j}}\right)_{\phi(p)} + \displaystyle\sum_{j=1}^{m}\frac{\partial b_{j}}{\partial y^{i}}(p)J_{2}\left(\frac{\partial}{\partial v^{j}}\right)_{\phi(p)}
\end{array}
\]
\begin{equation}\label{eq 04:eq 04}
= \displaystyle\sum_{j=1}^{m}\frac{\partial a_{j}}{\partial y^{i}}(p)\frac{\partial}{\partial v^{j}}|_{\phi(p)} - \displaystyle\sum_{j=1}^{m}\frac{\partial b_{j}}{\partial y^{i}}(p)\frac{\partial}{\partial u^{j}}|_{\phi(p)} 
\end{equation}

Now equation ({\ref{eq 02:eq 02}}) = equation ({\ref{eq 04:eq 04}}) if and only if $$\displaystyle\frac{\partial a_{j}}{\partial x^{i}} = \displaystyle\frac{\partial b_{j}}{\partial y^{i}},$$ that is, if and only if $$\displaystyle\frac{\partial u^{j}}{\partial x^{i}} = \frac{\partial v^{j}}{\partial y^{i}}$$\\

and equation ({\ref{eq 01:eq 01}}) = equation ({\ref{eq 03:eq 03}}) if and only if $$\displaystyle\frac{\partial a_{j}}{\partial y^{i}} = -\frac{\partial b_{j}}{\partial x^{i}},$$ that is, if and only if $$\displaystyle\frac{\partial v^{j}}{\partial x^{i}} = -\frac{\partial u^{j}}{\partial y^{i}}$$ which are Cauchy-Riemann equations. Thus $\phi$ is holomorphic if and only if $\phi_{*}\circ J_{1} = J_{2}\circ\phi_{*}$ as required. This completes the proof of the theorem. $\blacksquare$\\

\section{Almost complex manifolds}

If $(M , \omega)$ is a symplectic manifold, an almost complex structure $J$ on $M$ is said to be compatible if whenever $m\in M$ and $g_{m}:T_{m}M\times T_{m}M\rightarrow \mathbb{R} $, then $g_{m}(u , v):=\omega_{m}(u , Jv)$ defines a Riemannian metric on $M$, for all $u,v\in T_{m}M$.

\begin{proposition}\label {pr:pr01}
For every symplectic manifold $(M , \omega)$ there exists an almost complex structure $J$ and a Riemannian metric $g(\cdot , \cdot)$ on $M$ such that for each $m\in M$ we have $\omega_{m}(u , Jv) = g_{m}(u , v)$ for all $u, v\in T_{m}M$.

\end{proposition}
See \cite [Proposition 5] {H-Z94}.

\vspace{0.3cm}

Note that we can also write the compatibility condition in the form
\[
\begin{array}{cc}
\omega_{m}(u, v) = g_{m}(Ju , v), \hspace{0.4cm} u,v\in T_{m}M
\end{array}
\]

\begin{proposition}\label{pr:pr02}
Let $G$ be a compact Lie group and $\Phi:G\times M\rightarrow M$ a symplectic action of $G$ on the symplectic manifold $(M , \omega)$. Let $g(\cdot , \cdot)$ be an invariant metric on $M$ and $A$ a field of endomorphisms of $TM$, that is, $A:TM\rightarrow TM$ such that for each $m\in M$ we have $\omega_{m}(X , Y) = g_{m}(A_{m}X , Y)$, $X, Y\in T_{m}M$, then $A$ is $G$-invariant.
\end{proposition}

\proof  Let $a\in G$, $m\in M$. Suppose further that $X, Y$ are vectors such that $X\in T_{m}M$ $Y\in T_{\Phi_{a}(m)}M$. Then we have:

\[
\begin{array}{cll}
g_{\Phi_{a}(m)}(T_{m}\Phi_{a}\circ A_{m}X , Y)
&=& g_{m}(A_{m}X , (T_{m}\Phi_{a})^{-1}Y)\\
&=& \omega_{m}(X , (T_{m}\Phi_{a})^{-1}Y)\\
&=& \omega_{\Phi_{a}(m)}(T_{m}\Phi_{a}\cdot X , Y)\\
&=& g_{\Phi_{a}(m)}(A_{\Phi_{a}(m)}\circ (T_{m}\Phi_{a})X , Y)
\end{array}
\]

Thus $T_{m}\Phi_{a}\circ A_{m} = A_{\Phi_{a}(m)}\circ T_{m}\Phi_{a}$

This proves the proposition.

\begin{proposition}\label{pr:pr03}
Let $(M , \omega)$ be a symplectic manifold with a compatible almost complex structure $J$. If $G$ is a group of isometries of $M$ acting in a symplectic way, then the compatible almost complex structure $J$ is $G$-invariant.
\end{proposition}

\proof Let $g$ be a Riemannian metric on $M$ such that for each $x\in M$, we have $g_{x}(Ju,v) = \omega_{x}(u,v)$ for all $u,v\in T_{x}M$. Then, for all $x\in M$ we have:
\[
\begin{array}{cll}
 g_{x}(Ju,v) &=& \omega_{x}(u,v) = \Phi_{a}^{*}\omega_{x}(u,v)\\
 &=& \omega_{\Phi_{a}(x)}(T_{x}\Phi_{a}u,T_{x}\Phi_{a}v)\\
 &=& g_{\Phi_{a}(x)}(JT_{x}\Phi_{a}u,T_{x}\Phi_{a}v)\\
 &=& g_{\Phi_{a}^{-1}\circ\Phi_{a}(x)}(T_{x}\Phi_{a}^{-1}\circ J\circ T_{x}\Phi_{a}u,v)\\
 &=& g_{x}(T_{x}\Phi_{a}^{-1}\circ J\circ T_{x}\Phi_{a}u,v)
\end{array}
\]
for all $u,v\in T_{x}M$\\

Thus $Ju = T_{x}\Phi_{a}^{-1}\circ J \circ T_{x}\Phi_{a}u$ which gives $T_{x}\Phi_{a}\circ J = J\circ T_{x}\Phi_{a}$ $\blacksquare$

\vspace{0.5cm}
\subsection{Symplectic reduced space}

\begin{definition}
Let $(M,\omega)$ be a symplectic manifold and $G$ a Lie group. Let $\Phi:G\times M\rightarrow M$ be a Hamiltonian action of $G$ on $M$. Let $\mu:M\rightarrow\mathfrak{g}^{*}$ be the $Ad^{*}$-equivariant momentum mapping of the action and $\beta\in\mathfrak{g}^{*}$ a regular value of $\mu$. We define the symplectic reduced space of the $G$-action on $M$ to be
\[
\begin{array}{ccc}
M_{\beta}:= \mu^{-1}(\beta)/G_{\beta},
\end{array}
\]
where $G_{\beta}$ is the isotropy subgroup of $\beta$.

\end{definition}

\begin{description}
\item   (i) Since $\beta\in\mathfrak{g}^{*}$ is a regular value of $\mu$, the inverse image $\mu^{-1}(\beta)$ is a submanifold of $M$ of dimension $\dim M-\dim G$.
\item   (ii)   If the action of $G_{\beta}$ on $\mu^{-1}(\beta)$ is free and proper then the reduced space $M_{\beta} = \mu^{-1}(\beta)/G_{\beta}$ is a manifold of dimension $\dim M- \dim G-\dim G_{\beta}$
\end{description}
See \cite [p.~124]{Mar74}. In this case, the projection map $\pi_{\beta}:\mu^{-1}(\beta)\rightarrow \mu^{-1}(\beta)/G_{\beta}$ is a smooth submersion and $\omega_{\beta}$, where \begin{equation}\label{eq 05:eq 05}
\pi_{\beta}^{*}\omega_{\beta} = i_{\beta}^{*}\omega 
\end{equation} is the unique symplectic form on the reduced space $M_{\beta}$ by the Marsden-Weinstein-Meyer reduction theorem (See \cite [pp.~298-299] {Abr78}), 

with $i_{\beta}:\mu^{-1}(\beta)\rightarrow M$ the inclusion map and $\pi_{\beta}:\mu^{-1}(\beta)\rightarrow \mu^{-1}(\beta)/G_{\beta}$ the quotient map. That is, if $x$ is a point in $\mu^{-1}(\beta)$ so that $\pi_{\beta}(x) = [x]$ is a point on the quotient space $\mu^{-1}(\beta)/G_{\beta}$ and $u\in  T_{x}(\mu^{-1}(\beta))$ be a tangent vector so that $[u]\in T_{[x]}(\mu^{-1}(\beta)/G_{\beta})$ identified with quotient of tangent spaces $T_{x}(\mu^{-1}(\beta))/T_{x}(G_{\beta}\cdot x)$, then the equation ({\ref{eq 05:eq 05}}) is equivalent to the following
\[
\begin{array}{cc}
\omega_{\beta}([x])([u] , [v]) = \omega(x)(u , v)
\end{array}
\]

for all $u , v \in T_{x}(\mu^{-1}(\beta))$. (See \cite [p.~15] {Mar07}).\\

 Let $g_{M}$ be a Riemannian metric on the symplectic manifold $(M,\omega)$ and let $J_{M}$ be an almost complex structure such that $\omega(\cdot , \cdot) = g_{M}(J_{M}\cdot , \cdot)$, then for $u , v\in T_{x}(\mu^{-1}(\beta))$ we have
\[
\begin{array}{cll}
i^{*}\omega(x)(u , v)	&=& \omega(x)(i_{*}u , i_{*}v) \\
						&=& g_{M}(x)(J_{M}(i_{*}u) , i_{*}v)\\
						&=& g_{M}(x)(J_{M}u , v) \\
						&=& g_{M}(x)i_{*}(J_{M}u) , i_{*}v) = i^{*}g_{M}(x)(J_{M}u , v),

\end{array}
\]
one has $i^{*}\omega(\cdot,\cdot) = i^{*}g_{M}(J_{M}\cdot,\cdot)$.

\begin{theorem}\label{th:th01}
Let $(M,\omega)$ be a symplectic manifold and $G$ a Lie group of isometries of $M$ whose action on $M$ is a hamiltonian action. Let $\mu :M\rightarrow\mathfrak{g}^{*}$ be the $Ad^{*}$-equivariant momentum mapping of the action, where $\mathfrak{g}^{*}$ is the dual of the Lie algebra of $G$. Let $\beta\in\mathfrak{g}^{*}$ be a regular value of $\mu$ and $G_{\beta}$ the isotropy subgroup of $\beta$ which acts freely and properly on $\mu^{-1}(\beta)$. Then there exists a Riemannian metric $g_{\beta}$ on the reduced space $\mu^{-1}(\beta)/G_{\beta}$ such that the projection map $\pi_{\beta}:\mu^{-1}(\beta)\rightarrow\mu^{-1}(\beta)/G_{\beta}$ is a Riemannian submersion. That is, $\pi_{\beta}^{*}g_{\beta} = i^{*}g_{M}$ where $g_{M}$ is a Riemannian metric on $M$ and $i:\mu^{-1}(\beta)\rightarrow M$ is the inclusion map.
\end{theorem}

\proof

Let $\pi:\mu^{-1}(\beta)\rightarrow \mu^{-1}(\beta)/G_{\beta}$ be the projection onto the reduced space. For convenience we shall write $ M_{\beta}$ for $\mu^{-1}(\beta)$ and $B_{\beta}$ for $\mu^{-1}(\beta)/G_{\beta}$. If $x\in B_{\beta}$ then $\pi^{-1}(x)$ is called the fibre over $x$. If $m\in \pi^{-1}(x)$ then $\pi^{-1}(x) = \lbrace gm:g\in G ~~\textrm{and}~~ \pi(gm) = x\rbrace$ is the fibre through $m$. Since $G_{\beta}$ acts freely and properly on $\mu^{-1}(\beta)$, the projection $\pi:M_{\beta}\rightarrow B_{\beta}$ is a submersion. (See \cite [pp.~298-299] {Abr78}). But $\pi$ is constant on $\pi^{-1}(x)$, for each $x\in B_{\beta}$, that is, $\pi(\pi^{-1}(x))=\lbrace x\rbrace$, so, if $u\in T_{m}\pi^{-1}(x)$ for $m\in \pi^{-1}(x)$ then $d\pi_{m}(u)=0$. That is, $T_{m}\pi^{-1}(x) = \ker{d\pi_{m}} = V(M _{\beta})_{m}$ is the set of vertical vectors.  Let $H(M_{\beta})_{m}$ be the orthogonal complement of $V(M_{\beta})_{m}$, then $T_{m}M_{\beta}$ decomposes into a direct sum
\[
\begin{array}{ccc}
T_{m}M_{\beta} = H(M_{\beta})_{m}\oplus V(M_{\beta})_{m},
\end{array}
\]
with $H(M_{\beta})_{m}\cap V(M_{\beta})_{m} = \lbrace 0 \rbrace$. Thus, if $X\in T_{m}M_{\beta}$ then $X=Y+Z$ with $Y\in H(M_{\beta})_{m}$ and $Z\in V(M_{\beta})_{m}$. It follows that $d\pi_{m}(X) = d\pi_{m}(Y)$. So, if $X\not\in V(M_{\beta})_{m}$ then $d\pi_{m}(X)\neq 0$ and $d\pi_{m}(X)\in T_{[m]}B_{\beta}$, where $[m]=\pi(m)$. Thus, for each $X\in H(M_{\beta})_{m}$ we have $d\pi_{m}(X)\in T_{[m]}B_{\beta}$. Let $d\pi_{m}|_{H_{m}}$ be the restriction of $d\pi_{m}$ to $H(M_{\beta})_{m}$, the space of horizontal vectors. Since $\pi$ and $d\pi$ are surjective (\cite [p.~299]{Abr78}), then $d\pi_{m}|_{H_{m}}$ is surjective and it is linear. But $\ker{d\pi_{m}|_{H_{m}}} = \lbrace 0 \rbrace$, so if $[u]\in T_{[m]}B_{\beta}$, there is a unique $u\in H(M_{\beta})_{m}$ such that $d\pi_{m}|_{H_{m}}(u) = [u]$. That is the map $d\pi_{m}|_{H_{m}}$ is also injective. It follows therefore that the map
\[
\begin{array}{ccc}
d\pi_{m}|_{H_{m}}:H(M_{\beta})_{m}\rightarrow T_{[m]}B_{\beta}
\end{array}
\]
is an isomorphism of vector spaces. Because of this isomorphism we shall write the tangent vectors of $T_{[m]}B_{\beta}$ as say $w$ instead of $[w]$ when we refer to the restriction map $d\pi_{m}|_{H_{m}}$.

Now by proposition \ref{pr:pr01}, there is an almost complex structure $J_{M}$ on $M$ and a Riemannian metric $g_{M}$ on $M$ such that if $X,Y\in T_{m}M$, then $\omega(m)(X,J_{M}Y) = g_{M}(m)(X,Y)$.

Let $v , w\in T_{x}B_{\beta}$. Then there exists unique vectors $\tilde{v} , \tilde{w}\in H(M_{\beta})_{m}$, $m\in \pi^{-1}(x)$ such that $d\pi_{m}|_{H_{m}}(\tilde{v}) = v$ and $d\pi_{m}|_{H_{m}}(\tilde{w}) = w$.

Define a metric $h$ on $T_{x}B_{\beta}$ by $h_{x}(v , w) = i^{*}g_{M}(\tilde{v} , \tilde{w})$. We shall show that the assignment $x\mapsto h_{x}$ smoothly depends on $x$. First note that if $m_{1}, m_{2}\in \pi^{-1}(x)$ then there is an isometry $f\in G$ with $f(m_{1}) = m_{2}$ and $\pi\circ f = \pi$. (See \cite [proposition 2.20]{G-H-L87}). We then have $T_{f(m_{1})}\pi\circ T_{m_{1}}f = T_{m_{1}}\pi$. Thus, $T_{m_{1}}f$ is an isometry between $H(M_{\beta})_{m_{1}}$ and $H(M_{\beta})_{m_{2}}$. This shows that $h_{x}$ does not depend on the choice of $m$ in the fibre $\pi^{-1}(x)$.

Let $m\mapsto p_{m}$ be a smooth assignment of the orthogonal projection $p_{m}:T_{m}M_{\beta}\rightarrow H(M_{\beta})_{m}$ of $T_{m}M_{\beta}$ onto $H(M_{\beta})_{m}$. Since $d\pi_{m}|_{H_{m}}$ is an isomorphism, $\pi$ is a local diffeomorphism. Let $\sigma$ be the local section of $\pi$. If $U$ is an open subset of $B_{\beta}$ and $x\in U$, let $v',w'\in T_{\sigma(x)}M_{\beta}$, then $h_{x}(v , w) = i^{*}g_{M}(\sigma(x))(p_{\sigma(x)}v' , p_{\sigma(x)}w')$, where $p_{\sigma(x)}v' = \tilde{v}\in H(M_{\beta})_{\sigma(x)}$ and $p_{\sigma(x)}w' = \tilde{w}\in H(M_{\beta})_{\sigma(x)}$. As the right-hand side is the composition of smooth maps we conclude that $x\mapsto h_{x}$ is smooth and $d\pi_{m}|_{H_{m}}:H(M_{\beta})_{m}\rightarrow T_{\pi(m)}B_{\beta}$ is an isometry. By this construction we have shown that $\pi:\mu^{-1}(\beta)\rightarrow \mu^{-1}(\beta)/G_{\beta}$ is a Riemannian submersion.

\begin{definition}
An almost Hermitian manifold is an almost complex manifold $(M , J)$ with a chosen Riemannian structure $g_{M}$ such that $g_{M}(JX,JY)=g_{M}(X,Y)$ for all $X,Y\in TM$
\end{definition}

\begin{definition}
Let $(M,J_{M})$ and $(N,J_{N})$ be almost Hermitian manifolds, a map $\Phi:M\rightarrow N$ is called almost complex if it commutes with almost complex structures, that is, if $\Phi_{*}\circ J_{M} = J_{N}\circ\Phi_{*}$
\end{definition}

An almost complex mapping between almost Hermitian manifolds which is also a Riemannian submersion is called a almost Hermitian submersion.

\begin{proposition}
Let $\Phi:M\rightarrow N$ be an almost Hermitian submersion, then the horizontal and the vertical distributions determined by $\Phi$ are $J_{M}$-invariant. That is
\[
\begin{array}{ccc}
J_{M}\lbrace V(M)\rbrace = V(M)\\
J_{M}\lbrace H(M)\rbrace = H(M)
\end{array}
\]
\end{proposition}

\proof 
Let $(M,J_{M},g_{M})$ and $(N,J_{N},g_{N})$ be two almost Hermitian manifolds and $\Phi:M\rightarrow N$ an almost Hermitian submersion. Then $\Phi$ is an almost complex mapping and we have $\Phi_{*}\circ J_{M} = J_{N}\circ \Phi_{*}$. Let V be a vertical vector, then $\Phi_{*}V = 0$ since $V\in \ker\Phi_{*}$. We now have $\Phi_{*}(J_{M}V) = J_{N}(\Phi_{*}V = 0$. Thus $\Phi_{*}(J_{M}V) = 0$ which shows that $J_{M}V$ is a vertical vector. If now $X$ is a horizontal vector then for any vertical vector $V$ we have $g_{M}(X,V) = 0$ since they belong to orthogonal complement subspaces. We then have $g_{M}(J_{M}X,V) = g_{M}(J_{M}^{2}X,J_{M}V) = -g_{M}(X,J_{M}V) = 0$. Thus, $J_{M}X$ is horizontal vector.
See \cite [p.~151] {Wat76}.

\begin{definition}
Let $\Phi:M\rightarrow N$ be an almost Hermitian submersion. A horizontal vector field $X$ on $M$ is called a basic vector field if there is a smooth vector field denoted by $X_{*}$ on $N$ such that $X$ and $X_{*}$ are $\Phi$-related.
\end{definition}

We shall first state the difficulties that may arise with regard to the almost complex structure.\\

If for example a the symplectic manifold $(M,\omega)$ is a real manifold and $g_{M}$ is the Riemannian structure on $M$ such that $\omega(\cdot,\cdot)=g_{M}(J\cdot,\cdot)$, let $X$ be a vector field on $M$, then $0 = \omega(X,X) = g_{M}(JX,X)$. That is, $g_{M}(JX,X) = 0$ and since $g_{M}$ is positive definite we conclude that $JX$ is orthogonal to $X$. Thus, if $X$ is a horizontal vector field then $JX$ belong to the orthogonal complement which in this case is the vertical distribution. Therefore, even if $X$ is a basic vector field there is no guarantee that $JX$ will be a basic vector field.\\

Another difficulty arises from the push-forward of the almost complex structure. Even when the kernel of the differential of $\pi_{\beta}$ is preserved by $J$, there need not be an almost complex structure on the image $\pi_{\beta}(M_{\beta})$ which make $d\pi_{\beta}$ complex linear as the following example shows.

Consider the twistor fibration (see \cite {Bur86}).
\[
\begin{array}{cc}
\pi:\mathbb{C}P^{3}\rightarrow \mathbb{H}P^{1} = S^{4}\\
\mathbb{C}\cdot v\mapsto \mathbb{H}\cdot v, \hspace{0.5cm} v\in\mathbb{C}^{4}
\end{array}
\]
which sends a complex line through the origin in $\mathbb{C}^{4}$ to its quaternionic span in $\mathbb{H}^{2}$. For each point $x\in \mathbb{H}P^{1}$, the inverse image $\pi^{-1}(x)$ are complex lines in $\mathbb{C}P^{3}$. Thus the fibers of $\pi$ are holomorphic submanifolds of $\mathbb{C}P^{3}$ which are compact and connected. However, it has been proved that $\mathbb{H}P^{1}$ does not admit any almost complex structure. This shows that the push-forward of an almost complex structure by a submersion does not necessarily yield an almost complex structure on its image for which the differential of the map is complex linear. (See \cite [p.~8] {Bai-}) for the details of this example.

Another example of this phenomenon is found among covering maps of smooth manifolds $\pi:E\rightarrow B$ where $E$ has a complex structure. The immediate example is the covering map $\mathbb{C}P^{1}\rightarrow\mathbb{R}P^{2}$. It is immediate that $\mathbb{R}P^{2}$ does not admit any complex structure since it is not orientable.
\vspace{0.2cm}

We are now able to state and prove the main result of this work.

\begin{theorem}
Let $(M,\omega)$ be a symplectic manifold and $G$ a Lie group of isometries of $M$. Let $\Phi:G\times M\rightarrow M$ be a hamiltonian action of $G$ on $M$ with $Ad^{*}$-equivariant momentum mapping $\mu :M\rightarrow\mathfrak{g}^{*}$. Let $\beta\in\mathfrak{g}^{*}$ be a regular value of $\mu$ and $G_{\beta}$ be the isotropy subgroup of $\beta$ acting freely and properly on $\mu^{-1}(\beta)$. Given a compatible almost complex structure $J_{M}$ on $M$ and a Riemannian metric $g_{M}$ which satisfies the compatibility condition $\omega(X,Y) = g_{M}(J_{M}X,Y)$ for all $X,Y\in TM$, let $\omega_{\beta}$ be the reduced symplectic form on the reduced symplectic manifold $\mu^{-1}(\beta)/G_{\beta}$. Then there exists an almost complex structure $J_{\beta}$ and a Riemannian metric $g_{\beta}$ on the reduced space $\mu^{-1}(\beta)/G_{\beta}$ which make $\pi:\mu^{-1}(\beta)\rightarrow\mu^{-1}(\beta)/G_{\beta}$ a Riemannian submersion and satisfies the condition $\omega_{\beta}([u],[v]) = g_{\beta}(J_{\beta}[u],[v])$ for all $[u],[v]\in T(\mu^{-1}(\beta)/G_{\beta})$ if and only if $\pi:\mu^{-1}(\beta)\rightarrow \mu^{-1}(\beta)/G_{\beta}$ is an almost complex mapping.
\end{theorem}

\proof ({\bf of the main theorem}) Let $h_{\beta}$ be the Riemannan metric on $\mu^{-1}(\beta)/G_{\beta}$ as in theorem \ref{th:th01}. Since $\mu^{-1}(\beta)/G_{\beta}$ is a symplectic manifold, there is an almost complex structure $J_{\beta}$ and a Riemannian metric $g_{\beta}$ such that if $[u],[v]\in T(\mu^{-1}(\beta)/G_{\beta})$ then $\omega_{\beta}([u],[v]) = g_{\beta}(J_{\beta}[u],[v])$, see proposition \ref{pr:pr01}. It is sufficient to find a condition for which $h_{\beta} = g_{\beta}$. Let $x\in\mu^{-1}(\beta)/G_{\beta}$, we have seen from theorem \ref{th:th01} that if $m\in \pi^{-1}(x)$, then $\pi^{-1}(x)=\lbrace gm:g\in G\rbrace$ is the fibre through $m$ . The tangent space to the fibre $T_{m}(\pi^{-1}(x))$ is the kernel of the differential of $\pi$ at $m$. That is, $\ker{d\pi_{m}} = T_{m}(\pi^{-1}(x))$. We have classified this tangent space as the set of vertical vectors of the Riemannian submersion $\pi$. We also have by the Symplectic Reduction Theorem (see \cite [p.~15] {Mar07}) that $(T_{m}(\mu^{-1}(\beta)))^{\omega}=T_{m}(G\cdot m)$. But $G\cdot m = \lbrace gm:g\in G\rbrace = \pi^{-1}(x)$ is the fibre through $m$. So if $X\not\in T_{m}(\pi^{-1}(x))$ then there is a $Y\in T_{m}(\mu^{-1}(\beta))$ such that $\omega(X,Y)\neq 0$. That is,

\begin{equation}\label{eq 06:eq 06}
\omega(m)(X,Y)=\omega_{\beta}([m])([X],[Y]) = g_{\beta}([m])(J_{\beta}[X],[Y])\neq 0.
\end{equation}

But we also have that

\begin{equation}\label{eq 07:eq 07}
\omega(m)(X,Y)= g_{M}(m)(J_{M}X,Y) = h_{\beta}(\pi(m))(\pi_{*}(J_{M}X),\pi_{*}Y)
\end{equation}

by theorem \ref{th:th01}. In particular, if $X$ and $Y$ are basic vector fields then equation ({\ref{eq 06:eq 06}}) and ({\ref{eq 07:eq 07}}) imply that
\[
\begin{array}{cll}
g_{\beta}(\pi(m))(J_{\beta}(\pi_{*}X),\pi_{*}Y) &=& h_{\beta}(\pi(m))((J_{M}X)_{*},Y_{*})\circ\pi\\
												&=& \pi^{*}h_{\beta}(m)(JX,Y) \\
												&=& h_{\beta}(\pi(m)(\pi_{*}(J_{M}X),\pi_{*}Y).
\end{array}
\]
But this relation holds if and only if $J_{\beta}(\pi_{*}X) = \pi_{*}(J_{M}X)$, if and only if $\pi_{*}\circ J_{M} = J_{\beta}\circ\pi_{*}$, if and only if $\pi$ is an almost complex mapping.$\blacksquare$
\\
\\
{\bf Aknowledgement}: Augustin Batubenge is grateful to Francois Lalonde, along with the department of mathematics and Statistics of the Universit\'e de Montr\'eal for hosting him during the last two years as a visiting researcher. \\
Wallace Haziyu is thankful to Doctors I.D. Tembo, M. Lombe, and A. Ngwengwe for encouragements they rendered towards his doctoral research.

%\newpage

\vspace{0.2cm}

{\bf AUTHORS}:\\

{\bf Augustin T. Batubenge*}\\
{\small (corresponding author)}\\
Department of Mathematics and Statistics\\
University of Montreal and\\
University of Zambia\\
Email: a.batubenge@gmail.com\\

\vspace{0.2cm}

{\bf Wallace Haziyu}\\
Department of Mathematics and Statistics\\
University of Zambia\\
Email: whaziyu@unza.zm


\begin{thebibliography}{99}

\bibitem{Abr78}
R. Abraham and J. E. Marsden.
 \emph{\textit{Foundations of
Mechanics, (Second Edition)}}. The Benjamin/Cummings Publishing
Company, Inc, 1978.

\bibitem{Aud04}
M. Audin. \emph{Torus Action On Symplectic Manifolds}.~{Second Revised Version}. Birkhauser Verlag, Berlin, 2004.

\bibitem{Bai-}
Paul Baird. \emph{An Introduction To Twistors}

\bibitem{Bur86}
H.E. Burstall. \emph{\textrm{A twistor Discription of Harmonic Maps of a 2-Sphere into a Grassman}}. Math.(ematiche) Ann.(alem) 274~ Springer-Verlag, New York, 1986.
\bibitem{Che05}
Bang-Yen Chen.
\emph{Examples and classification of Riemannian submersions satisfying a basic equality} in Bulletin of the Australian Mathematical Society. V.72 Number 3, (2005), 391-402.

\bibitem{G-H-L87}
S. Gallot, D. Hulin, J. Lafontaine. \emph{Riemannian Geometry}. Springer-Verlag, Berlin Heidelberg, 1987.
\bibitem{H-Z94}
H. Hofer, E. Zehnder. \emph{Symplectic Invariants and Hamiltonian Dynamics}. Birkhauser Verlag, Switzerland, 1994.
\bibitem{Mar07}
J.E.Marsden,G.Misiolek, J.P. Ortega, M.Perlmutter, T.S.Ratiu. \emph{Hamiltonian Reduction by Stages}. Springer-Verlag, Berlin Heidelberg, 2007.

\bibitem{Mar74}
J. Marsden and A. Weinstein.
\emph{Reduction of Symplectic Manifolds with Symmetry} In Reports on Mathematical Physics Vol 15 (1974) Number 1, 121 - 129.

\bibitem{Mat72}
Y. Matsushima. \emph{Differentiable Manifolds}. Marcel Dekker Inc, New York, 1972.
\bibitem{Bar66}
B. O'Neil. \emph{The Fundamental Equations of a Submersion}. Michigan Math. J13 (1966), 459-469.
\bibitem{Wat76}
B. Watson. \emph{ Almost Hermitian Submersions}. Journal of Differential Geometry II (1976), 147-165.


\end{thebibliography}
\end{document}